\pdfoutput=1

\documentclass[preprint,review,12pt]{elsarticle}

\usepackage{graphicx}

\usepackage{amssymb}
\usepackage{amsmath}
\usepackage{amsthm}
\usepackage{url}
\usepackage{slashbox}

\usepackage{dsfont}

\newtheorem{theorem}{Theorem}

\newtheorem{definition}{Definition}

\journal{arXiv}

\begin{document}

\begin{frontmatter}



\title{Upper-critical graphs (complete $k$-partite graphs)}


\author{Jos\'e Antonio Mart\'in H.\corref{cor1}} \ead{jamartinh@fdi.ucm.es}
\cortext[cor1]{address: Facultad de Inform\'atica, Universidad Complutense de Madrid,
C. Prof. Jos\'e Garc\'ia Santesmases, s/n., 28040, Madrid, (Spain). Tel.: +34 91.394.7656, Fax: +34  91.394.7510}

\address{Faculty of Computer Science, Complutense University of Madrid, Spain}

\begin{abstract}
This work introduces the concept of \emph{upper-critical graphs}, in a complementary way of the conventional (lower)critical graphs: an element $x$ of a graph $G$ is called \emph{critical} if $\chi(G-x)<\chi(G)$. It is said that $G$ is a \emph{critical graph} if every element (vertex or edge) of $G$ is critical. Analogously, a graph $G$ is called \emph{upper-critical} if there is no edge that can be added to $G$ such that $G$ preserves its chromatic number, i.e. \{ $e \in E(\overline{G}) \; | \; \chi(G+e) = \chi(G)$ \} $=$ $\emptyset$. We show that the class of upper-critical graphs is the same as the class of complete $k$-partite graphs. A characterization in terms of hereditary properties under some transformations, e.g. subgraphs and minors and in terms of construction and counting is given.
\end{abstract}



\end{frontmatter}

\section{Preliminary definitions and basic terminology}
\label{sec:1}

Unless we state it otherwise, all graphs in this work are connected and simple (finite, and have
no loops or parallel edges).

Partitioning the set of vertices $V(G)$ of a graph $G$ into separate classes, in such a way that no two adjacent vertices are
grouped into the same class, is called the vertex coloring problem. In order to distinguish such classes, a set of colors $C$
is used, and the division into these \emph{(color)-classes} is given by a coloring $c: V(G)\rightarrow \{1...k\}$, where $c(x) \neq c(y)$ for all $xy's$ belonging to the set of edges $E(G)$ of $G$, where $c^*(x)$ is the color-class of $x$ (i.e. the color-class that contains vertex $x$).

Given a graph coloring over a graph $G$ with a set of colors $C$, if $C$ has cardinality $k$, then $c$ is a \emph{k-coloring} of $G$. The \emph{chromatic number} of a graph $\chi(G)$ is the minimum number of different colors which are necessary for coloring the vertices of a graph $G$. If $\chi(G) \leq k$ then $G$ is \emph{k-colorable} (i.e. can be colored with $k$ different colors) and if $\chi(G)= k$ then $G$ is \emph{k-chromatic}. A graph is \emph{complete} if every two distinct vertices in the graph are adjacent. The complete graph of order $k$ is denoted by $K_k$.

The \emph{join} $ G_1 + G_2 $ of two different (vertex-disjoint) graphs, $G_1$ and $G_2$, is the graph that has vertex set $V (G_1 + G_2) = V(G_1)\cup V(G_2)$ and edge set $E(G_1 + G_2) = E(G_1) \cup E(G_2) \cup \{uv | u \in V(G_1), v \in V(G_2)\}$, that is, each vertex of $G_1$ joined to each vertex of $G_2$.

An \emph{independent set} (also called \emph{stable set}) $I$ is a set of vertices of $G$ such that there are no edges between any two vertices in $I$. The \emph{complement} of $G$ is denoted by $\overline{G}$, e.g. An independent set $I$ of $G$ is a subgraph of $\overline{G}$ which is complete (i.e. a \emph{complete subgraph} or \emph{clique} of $\overline{G})$, therefore, we can denote an independent set of cardinality $n$ as $\overline{K}_n$.

A \emph{complete k-partite graph} $K_{n_1, n_2,...,n_k}$ is a graph isomorphic to the join of $k$ independent sets of cardinality $n_1, n_2,...,n_k$ respectively, i.e., $K_{n_1, n_2,...,n_k} = \overline{K}_{n_1} + \overline{K}_{n_2} + ... + \overline{K}_{n_k}$, e.g. $K_{2,2}$ is a square.

The set of all adjacent vertices to a vertex $x\in V(G)$ is called its \emph{neighborhood} and is denoted by $N_G(x)$. When it is clear to which graph we are referring to, we will use simply $N(x)$, omitting the graph. The \emph{closed neighborhood} of a vertex $x$, denoted by $N[G]$, includes also the vertex $x$, i.e. $N[x]=N(x)\cup \{x\}$. If $N(x)=N(y)$ we say that vertex $y$ is a copy of $x$ and vice versa.

The degree of a vertex $x$, $deg(x)$, is equal to the cardinality of its neighborhood $deg(x)=|N(x)|$. A \emph{complete vertex} is any $x \in V(G)$ such that $N[x] = V(G)$ and hence if every vertex is complete the graph is also complete $K_k$.

An \emph{edge contraction} denoted by $G/xy$ or $G/e$ is the process of replacing two adjacent vertices $x,y$ of $G$, i.e $xy \in E(G)$, by a new vertex $z$ such that $N(z) = N(x)\cup N(y)$.

A graph $H$ is called a \emph{minor} of the graph $G$ ($H\prec G$) if $H$ is isomorphic to a graph that can be obtained from a subgraph of $G$ by zero or more edge deletions, edge contractions or vertex deletions on a subgraph of $G$. In particular, $G$ is minor of itself.

A \emph{vertex identification} denoted by $G/x,y$ is the process of replacing two non-adjacent vertices $x,y$ of $G$, i.e $xy\notin E(G)$, by a new vertex $z$ such that $N(z) = N(x)\cup N(y)$.

A graph $H = G/S$ is called a \emph{contraction} of the graph $G$ if $H$ is isomorphic to a graph that can be obtained by zero or more vertex
identifications on $G$, where $S$ is the set of vertex identifications to obtain $H$. In particular, $G$ is a contraction of itself.

When $S$ induces a contraction $H$ of $G$ such that $H=K_k$, $H$ is called a \emph{collapse} of $G$. Thus, since a collapse $H=K_k$ is also a partition of all the vertices of $G$ in $k$ independent sets, $\chi(G)$ is the size of the minimum collapse of $G$ and if $k=\chi(G)$ then $H = K_k$ induces a $k$-coloring of $G$.

\emph{Critical graphs} were first studied by Dirac~\cite{dirac52,dirac52c,dirac53}. An element $x$ of a graph $G$ is called \emph{critical} if $\chi(G-x)<\chi(G)$. If all the vertices of a graph $G$ are critical we say that $G$ is \emph{vertex-critical} and if every element (vertex or edge) of $G$ is critical we say that $G$ is a \emph{critical graph}, specifically, if $\chi(G)=k$, $G$ is $k$-critical.

Examples of critical graphs in general are the \emph{complete graphs} $K_k$ of size $|V(G)|=k$, odd cycles are the only 3-critical graphs and odd wheels are just one case of 4-critical graphs. If $x$ is a critical vertex of a graph $G$ then $x$ is a color-class itself, that is, there is at least one $k$-coloring of $G$ where $x$ is the only vertex with the $c(x)$ color.

Also $k$-critical graphs possess the next well known properties:
\begin{itemize}
    \item G has only one component.
    \item G is finite
    \item Every vertex is adjacent to at least $k$-1 others.
    \item If G is ($k$-1)-regular, meaning every vertex is adjacent to exactly
          $k$-1 others, then G is either $K_k$ or an odd cycle.
    \item $|V(G)|\neq k+1$
    \item If G is different from $K_k$ then $|V(G)| \geq k+2$.
\end{itemize}

\section{Upper-critical graphs}
This work introduces the concept of \emph{upper-critical graphs}, in a complementary way of the conventional (lower)critical graphs.
Also, we show that a $k$-chromatic graph $G$ is upper-critical if and only if it is isomorphic to a complete $k$-partite graph, showing that the class of the upper-critical graphs coincides with the class of the complete $k$-partite graphs.

\begin{definition} A graph $G$ is called \emph{upper-critical} if there is no edge ($e$) that can be added to $G$ such that $G$ preserves its chromatic number, i.e.
$$\{e \in E(\overline{G}) \; | \; \chi(G+e) = \chi(G)\} \;=\; \emptyset $$
or alternatively,
$$\chi(G+e) > \chi(G) \;\; \forall e \in E(\overline{G}), \; \text{ or } \;\; G=K_k,$$
If $\chi(G)=k$, $G$ is called a $k$-chromatic upper-critical graph.
\end{definition}

\subsection{General structural properties}

\begin{theorem}\label{uc:uniquely} If $G$ is a $k$-chromatic upper-critical graph then every $k$-coloring of $G$ induces the same partition of vertices of $G$, in $k$ different color classes, i.e. $G$ is uniquely colorable (uniquely $k$-colorable, to be more precise).
\begin{proof}
Given a $k$-chromatic upper-critical graph $G$, i.e. $\chi(G+e) > \chi(G)$. Let us suppose that there are at least two different colorings of $G$ that induce two different partitions of $G$. So, there is at least a vertex $x$ that can be assigned to either a color-class $c_1$ or a color-class $c_2$. Hence, it is possible to add a new edge $e=xy$ to $G$ ($G+e$), from vertex $x$ towards an element of $c_2$ (say vertex $y$). But, in this case, $G+e$ has at least one $k$-coloring which is a contradiction since $G$ is upper-critical.
\end{proof}
\end{theorem}

\begin{theorem}\label{uc:nx} If $G$ is an upper-critical graph, $x$ a vertex of $G$ and $c$ a $k$-coloring of $G$ then:
$$N(x)= V(G) - \{y\in V(G) \;|\; c(x)\neq c(y) \}$$

\begin{proof} Suppose $N(x) \neq V(G) - \{y\in V(G) \;|\; c(x)\neq c(y) \}$. Then either $c$ is not a $k$-coloring of $G$ or there is a vertex $y\in V(G)$ such that $y \notin N(x)$ and $c(x)\neq c(y)$. Hence $\chi(G+xy) = \chi(G)$, which is a contradiction since $G$ is upper-critical.
\end{proof}
\end{theorem}

\begin{theorem}\label{uc:completesubgraph} If $G$ is a $k$-chromatic upper-critical graph then:
$$K_k \subseteq G$$
\begin{proof}
Since $G$ has $k$ color classes, let $x_1,x_2,x_3,...,x_k$ be $k$ vertices such that each $x_i$ belongs to a different color class $c_i$.
Now, the induced subgraph $x_1,x_2,x_3,...,x_k$ is a complete graph $K_k$ since, by theorem~\ref{uc:nx}, $x_ix_j \in E(G) \; \forall i \neq j$.
\end{proof}
\end{theorem}

\subsection{Hereditary properties of transformations, subgraphs and minors}

In this section we will proof that upper-critical graphs are closed under vertex deletion, vertex identification and edge contraction. And also, for particular cases, under vertex addition, edge addition and edge deletion, i.e., that given an upper-critical graph $G$, the next graphs are also upper-critical:
\begin{enumerate}
  \item $G-x$.
  \item $G/x,y$.
  \item $G/xy$.
  \item $G+x$; if $x$ is a copy of some vertex of $G$ or $x$ is a complete vertex.
  \item $G+e$; for $e=xy\notin E(G)$ if either:
  \begin{enumerate}
    \item $N[x]\neq K_k$ or,
    \item $|E(\overline{G})| \leq |V(G)|-\chi(G)$
  \end{enumerate}
  \item $K_k-\{e_1,...,e_{k-2} \}$; for some particular sequence of critical edges.
\end{enumerate}


Subsequent sections, will study how are the graphs obtained from $G$ by applying each particular graph transformation, e.g., decrease, increase or maintain its chromatic number or its order.

\begin{theorem}\label{uc:vertext-copy} If $G$ is an upper-critical graph and $xy \notin E(G)$ then:
$$ N(x) = N(y), \text{ i.e. vertex $y$ is a ``copy'' of vertex $x$} $$
\begin{proof}
Let us suppose there is a vertex $z$ in $N(y)$ but not in $N(x)$, i.e. $ N(x) \neq N(y)$, then:
\begin{itemize}
  \item $\chi(G+xz) > \chi(G)$ hence $c(x)=c(z)$ for every $\chi$-coloring $c$ of $G$.
  \item $\chi(G+xy) > \chi(G)$ hence $c(x)=c(y)$ for every $\chi$-coloring $c$ of $G$.
\end{itemize}
Hence, $c(z) = c(y)$ for every $\chi$-coloring $c$ of $G$, which is a contradiction since $yz \in E(G)$.
\end{proof}
\end{theorem}

\begin{theorem}\label{uc:copyplus} If $G$ is an upper-critical graph and $x$ a vertex of $G$ then:
$$ G + y \;\;\text{ is also upper-critical if $y$ is a copy of $x$} $$
\begin{proof} The case were $G$ is a complete graph is trivial. Let $G$ be an upper-critical graph different from a complete graph. Now, let $G+y$ be the graph obtained by adding a copy ($y$) of the vertex $x$ to $G$. Then, since $c(y)=c(x)$ in every $k$-coloring of $G$ then $G+y$ is $k$-chromatic and, since $\chi(G+y+e) > \chi(G) \text{ for all } e \notin E(G+y)$ then $G+y$ is upper-critical.
\end{proof}
\end{theorem}

\begin{theorem}\label{uc:minus} If $G$ is an upper-critical graph and $x$ a vertex of $G$ then:
$$ G-x \;\;\text{ is also upper-critical} $$
\begin{proof}
Let us suppose that $G$ is a $k$-chromatic upper-critical graph and $e$ a non-existent edge in $G$. Assume $G-x$ is not upper-critical. We can consider just two cases for $G-x$:

\begin{enumerate}
  \item If $\chi(G-x) = k-1$ then $(G-x)+e$ has a ($k$-1)-coloring (since $G-x$ is not upper-critical), hence there is a $k$-coloring of $(G-x)+e+x$, i.e. a $k$-coloring of $G+e$, which is a contradiction since $G$ is a $k$-chromatic upper-critical graph.

  \item If $\chi(G-x) = k  $ then $x$ is not a complete vertex, meaning that there is another vertex $y$ such that $N(x)=N(y)$ and $xy\notin E(G)$. Now, since $G-x$ is not upper-critical there are two vertices $u,v\in V(G-x)$ such that $uv\notin E(G-x)$ and there is a $k$-coloring ($c$) of $G-x$ where $c(u) \neq c(v)$, i.e., a $k$-coloring of $(G-x)+e$. However, since $N(x)=N(y)$ we can restore back vertex $x$ assigning to it the color $c(y)$ obtaining a $k$-coloring of $G$ such that there are two vertices $u,v\in V(G)$,  $uv\notin E(G)$ and $c(u)\neq c(v)$, i.e., we can obtain a $k$-coloring of $G+e$, which is a contradiction since $G$ is upper-critical.
\end{enumerate}

\end{proof}
\end{theorem}

\begin{theorem}\label{t5} If $G$ is an upper-critical graph and $x,y$ are two vertices of $G$ such that $xy \notin E(G)$ then:
$$G/x,y  \;\;\text{ is also upper-critical}$$

\begin{proof}
Since $G$ is upper-critical and $xy \notin E(G)$ then, by theorem~\ref{uc:vertext-copy}, vertex $y$ is a copy of $x$ and hence $G/x,y = G-x$. Now, since $G-x$ is upper-critical then $G/x,y$ is upper-critical.
\end{proof}

\end{theorem}

\begin{theorem}\label{t6} If $G$ is an upper-critical graph and $xy \in E(G)$ then:
$$G/xy  \;\;\text{ is also upper-critical}$$
\begin{proof} by induction on the number $n=V(G)$ of vertices of $G$.Let $e$ be some edge of $G$.

Base $n=2$: The $K_2$ graph is upper-critical by definition and $K_2/e = K_1$ and $K_1$ is upper-critical by definition.

Assume true for every upper-critical graph in $n$ vertices.

Proof for $n+1$: Let $G$ be an upper-critical graph on $n+1$ vertices such that $G/e$ is not upper-critical.

\begin{enumerate}
  \item If $G=K_k$ then $G/e$ is also a complete graph hence $G/e$ is upper-critical which is a contradiction.
  \item If $G\neq K_k$ then there are two cases:
  \begin{enumerate}
  \item $x$ is a complete vertex (this case includes all three cases: x, y or both). It is easy to see that $G/xy=G-x$, but $G-x$ is upper-critical, thus there is a contradiction.
  \item Since $G \neq K_k$ and $x$ is not a complete vertex: we can delete a copy ($z$) of $x$ from $G$, obtaining a new graph $H=G-z$.
  Since $G$ is upper-critical then $H$ is upper-critical and thus $(G-z)/e$ is upper-critical by the inductive hypothesis. But now since $z$ is a copy of $x$ then: $$(G-z)/e + z = G/e $$ is upper-critical by theorem~\ref{uc:copyplus} which is a contradiction.
\end{enumerate}
\end{enumerate}
\end{proof}
\end{theorem}

\begin{theorem}\label{uc:plusedge} If $G$ is an upper-critical graph, $x,y$ are two vertices of $G$ such that $xy \notin E(G)$ then:
$$G+xy  \;\;\text{ is also upper-critical if:}$$
\begin{enumerate}
  \item $N[x]\neq K_k$
  \item $|E(\overline{G})| \leq |V(G)|-k$
\end{enumerate}

\begin{proof} Since $xy \notin E(G)$ there is always a vertex $z$, adjacent to both $x$ and $y$.
\begin{enumerate}
  \item There is at least one $z$ which is not a complete vertex: Let $\dot{z}$ be a copy of $z$. Then $G+\dot{z}$ is upper-critical and so $(G+\dot{z})/x\dot{z}$. Now: $$(G+\dot{z})/x\dot{z} = G+xy = G+e$$
  \item Every $z$ is a complete vertex: Fix $G$ to be $k$-chromatic. Then $G+e$ has $K_{k+1}$ as a subgraph and $G+e+e_2$ will have $K_{k+2}$ as a subgraph, for a new edge $e_2$, so $\chi(G+e+e_2)>\chi(G+e)$. Hence $G+e$ is upper-critical.

\end{enumerate}
Therefore $G+e$ is upper-critical.
\end{proof}
\end{theorem}

\begin{theorem}\label{uc:minusedge} If $G=K_k$ then there is always a particular sequence (including the empty sequence) of critical edges such that: $$G=K_k-\{e_1,...,e_{k-2} \} \;\;\mbox{is also upper-critical}$$

\begin{proof} Let $G$ be a $2$-chromatic (bipartite) upper-critical graph on $N$ vertices then, by definition, $G+e$ is $(k+1)$-chromatic and, by theorem~\ref{uc:plusedge}, $G+e$ is also upper-critical, hence it is possible to add edges obtaining successive upper-critical graphs up to obtain a complete graph. Hence, given a particular complete graph there is always a particular sequence (including the empty sequence) of critical edges for obtaining successive upper-critical graphs with lower chromatic number up to a bipartite graph.
\end{proof}
\end{theorem}

\subsection{Construction, characterization and counting}

A procedure to obtain an upper-critical graph $G$, starting from any other graph, is to find a $k$-coloring of some $k$-chromatic graph $H$ and add the edges $x_ix_j$ to $H$ whenever the color class of vertex $x_i$ is different from the color class of $x_j$, i.e. $x_ix_j \in E(G) \; \forall i,j: c(x_i)\neq c(x_j)$. Thus, every coloring of any particular graph corresponds to a particular upper-critical graph.

Therefore, the total number of different colorings (excluding repeated partitions) of all graphs is equal to the the total number of upper-critical graphs.

Furthermore, there is a direct way to obtain any arbitrary upper-critical graph. Upper-critical graphs, contrary to (lower)critical-graphs, has a very easy general characterization and description, for instance,  the $n$-vertex list notation where $G=\{n_1,n_2,...,n_k\}$ is the $k$-chromatic upper-critical graph where the positive integers $\{n_1,n_2,...,n_k\}$ indicate de number of vertices in each color class ($c_i$) respectively.

It is immediate to see that the $n$-vertex list is unique in the sense that two upper-critical graphs share the same $n$-vertex list if and only if they are isomorphic.

From this it follows that it is possible to specify directly an arbitrary upper-critical graph using the $n$-vertex list notation subject to just one constraint:
\begin{equation}
|V(G)| = \sum_{i=1}^{k} n_i
\end{equation}
where $G$ is a $k$-chromatic upper-critical graph and $n_i$ is the number of vertices belonging to the $k$th color class.

As we can see, the definition of an upper-critical graph is the same as the definition of a complete $k$-partite graph. Theorem~\ref{kpartite:identity} formalizes this equivalence.
\begin{theorem}\label{kpartite:identity}
A $k$-chromatic graph $G$ is upper-critical if and only if it is isomorphic to a complete $k$-partite graph $H= K_{n_1, n_2,...,n_k}$.
\begin{proof}
It is clear that $H$ is $k$-chromatic and $$\{e \in E(\overline{H}) \; | \; \chi(H+e) = k \} \;=\; \emptyset, $$ so $H$ is an upper-critical graph. Also, since $G$ is $k$-chromatic it can be divided, by means of a $k$-coloring, into $k$-color classes (independent sets), which, by theorem~\ref{uc:uniquely}, induces always the same partition. Now, by theorem~\ref{uc:nx}, for every $x\in V(G)$ $N(x)= V(G) - \{y\in V(G) \;|\; c(x)\neq c(y) \}$ meaning that each element of an independent set of $G$ is completely joined to any other element of a different independent set of $G$. Therefore $G$ can be described as the join of $k$ independent sets, $G= K_{n_1, n_2,...,n_k}$, i.e.,  $G$ is a complete $k$-partite graph.
\end{proof}
\end{theorem}

Therefore, the above $n$-vertex list notation $G=\{n_1,n_2,...,n_k\}$ can be replaced by the standard notation $K_{n_1, n_2,...,n_k}$.

Table~\ref{tab:ucgraphs} shows the space of upper-critical graphs in  $(|V|,\chi)$ coordinates using the $n$-vertex list notation. A sample is shown up to $|V|=5$ and $\chi=5$.

\begin{table}[tb]
\caption{Upper-critical graphs on $|V(G)|$ vertices vs. Chromatic number $k=\chi(G)$}\label{tab:ucgraphs}
\centering
\begin{tabular}{|c|l|l|l|l|l|l|}
  \hline
  $V$ $\backslash$  $k$ & $1$ & $2$ & $3$ & $4$ & $5$ & ...\\
  \hline
  1 & $K_1$  &  &  &  & & \\
  \hline
  2 & &$K_2$ &  &  &  &\\
  \hline
  3 & &$\{1,2\}$ & $K_3$ &  &  &\\
  \hline
  4 & &$\{1,3\}$, $\{2,2\}$ & $\{1,1,2\}$ & $K_4$ &   &\\
  \hline
  5 & &$\{1,4\}$, $\{2,3\}$ & $\{1,1,3\}$, $\{1,2,2\}$ & $\{1,1,1,2\}$  & $K_5$ &\\
  \hline
  $\vdots$ & & $\vdots$ & $\vdots$ & $\vdots$  & $\vdots$ & $\ddots$ \\
  \hline
\end{tabular}
\end{table}

Furthermore, it is possible to count the number $P(N,k)$ of upper-critical graphs (i.e. complete $k$-partite graphs) given the number ($N=|V|$) of vertices and the chromatic number $(k)$. A \emph{partition} of a positive integer $N$ is a way to express $N$ as the sum of positive integers~\cite{Weisstein2010}, i.e. $4=2+2$ or $4=3+1$. Therefore, every partition of a positive integer defines a unequivocal upper-critical graph and viceversa.

Let $P(N,k)$ denote the number of ways of writing $N$ as a sum of exactly $k$ terms, then $P(N,k)$ can be computed from the next recurrence relation:
\begin{equation}
P(N,k)=P(N-1,k-1)+P(N-k,k),
\end{equation}

Let $\mathds{S}= |V| \times \chi$ be the space of all upper-critical graphs of order $N=|V|$ and chromatic number $k=\chi$, it is possible to go from a particular point $\mathds{S}(N_i,k_i)$ to another point $\mathds{S}(N_j,k_j)$ in $\mathds{S}$, i.e., traveling across table~\ref{tab:ucgraphs} cells. Note that to each point of $\mathds{S}$ could be associated several different upper-critical graphs sharing the same order and chromatic number. Theorem~\ref{uc:minusedge} formalizes this fact.

\begin{theorem}\label{thm:tabletravel} If $G$ is an upper-critical graph with $N=|V(G)|$ vertices and chromatic number $k$ then:
\begin{enumerate}
  \item $G+x$  $\in$
  \begin{enumerate}
    \item $\mathds{S}(N+1,k)$ (if $x$ is a copy of some non-complete vertex of $G$)
    \item $\mathds{S}(N+1,k+1)$ (if $x$ is a complete vertex)
  \end{enumerate}
  \item $G-x$    $\in$
  \begin{enumerate}
    \item $\mathds{S}(N-1,k)$ (if $x$ is a non-complete vertex)
    \item $\mathds{S}(N-1,k-1)$ (otherwise)
  \end{enumerate}
  \item $G/x,y$   $\in$ $\mathds{S}(N-1,k)$
  \item $G/xy$   $\in$
  \begin{enumerate}
    \item $\mathds{S}(N-1,k-1)$ (if $x$ and $y$ are both complete vertices)
    \item $\mathds{S}(N-1,k)$   (if $x$ or $y$ are non-complete vertices)
    \item $\mathds{S}(N-1,k+1)$ (if $x$ and $y$ are non-complete vertices)
  \end{enumerate} 
   \item $G+xy$  $\in$ $\mathds{S}(N,k+1)$ if either:
  \begin{enumerate}
    \item $N[x]\neq K_k$ or,
    \item $|E(\overline{G})| \leq |V(G)|-\chi(G)$
  \end{enumerate}
  \item $K_k-\{e_1,...,e_{m} \}$    $\in$  $\mathds{S}(N,k-m)$
\end{enumerate}
\begin{proof}$ $\\
\begin{enumerate}
  \item ...\\
       \begin{enumerate}
     \item if $x$ is a copy of some non-complete vertex of $G$ then $\chi(G)=\chi(G+x)$ since $x$ can be colored with the color of another vertex of $G$ and $G+x$ has $N+1$ vertices.
    \item Any graph plus a complete vertex increase both in vertex number and chromatic number, by one.
  \end{enumerate}
  \item ...\\
   \begin{enumerate}
    \item Since $x$ is a non-complete vertex then it is a copy of another vertex $y$, so $\chi(G) =\chi(G-x)$.
    \item Any graph minus a complete vertex decrease both in vertex number and chromatic number, by one.
    \\
  \end{enumerate}
  \item Since $x$ and $y$ are non-complete vertices then $G/x,y = G-x$.
  \item ...\\
   \begin{enumerate}
    \item if $x$ and $y$ are both complete vertices then $G/xy = G-x$
    \item if $x$ or $y$ are non-complete vertices then $G/xy = G-x$
    \item if $x$ and $y$ are non-complete vertices then $G/xy$ is equivalent to remove a copy of some vertex, $G-z$, and adding some edges, $G-z+e_1,e_2...$, which by definition will increase the chromatic number.
        \\
  \end{enumerate}
  \item By definition and theorem~\ref{uc:plusedge}.
  \item By definition and theorem~\ref{uc:minusedge}.
\end{enumerate}
\end{proof}
\end{theorem}


\bibliographystyle{plain}
\bibliography{gct}

\end{document}